\documentclass[graybox]{svmult}

\usepackage{type1cm}        % activate if the above 3 fonts are
\usepackage{makeidx}         % allows index generation
\usepackage{graphicx}        % standard LaTeX graphics tool
\usepackage{multicol}        % used for the two-column index
\usepackage[bottom]{footmisc}% places footnotes at page bottom

\usepackage{newtxtext}       % 
\usepackage[varvw]{newtxmath}       % selects Times Roman as basic font

\usepackage{lineno,hyperref}
\modulolinenumbers[5]

\usepackage{multirow}
\usepackage{slashbox}
\usepackage{subfigure}
\newcommand{\R}{\mathbb{R}}

\newcommand{\mb}[1]{\mathbf{#1}}

\newcommand\mbF{\mb{F}}
\newcommand\mbG{\mb{G}}
\renewcommand\div[1]{{\rm div\,} {#1}}

\newcommand\e[1]{{\rm e}^{#1}}
\renewcommand{\E}{\mathbb{E}}
\newcommand \Esp[1]{\E^Q\left[#1\right]}

\newcommand \dx[1]{{\rm d}#1}

\makeindex             % used for the subject index
\begin{document}

\title*{Second order finite volume IMEX Runge-Kutta schemes for two dimensional parabolic PDEs in finance}
\titlerunning{Second order FV IMEX RK for 2d PDEs in option pricing}
\author{J. G. L\'opez-Salas 
        and M. Su\'arez-Taboada
        and M. J. Castro
        and A. M. Ferreiro-Ferreiro
        and J. A. Garc\'ia-Rodr\'iguez
        }
\authorrunning{JG L\'opez-Salas, M Su\'arez-Taboada, A Ferreiro-Ferreiro, JA Garc\'ia, MJ Castro}

% \author[auth1]{J. G. L\'opez-Salas}\ead{jose.lsalas@udc.es}
% \address[auth1]{Department of Mathematics, Faculty of Informatics and CITIC, Campus Elvi\~na s/n, 15071-A Coru\~na (Spain)}

% \author[auth1]{M. Su\'arez-Taboada}\ead{maria.suarez3@udc.es}

% \author[auth2]{M. J. Castro\corref{cor1}}\ead{mjcastro@uma.es}
% \cortext[cor1]{Corresponding author.}
% \address[auth2]{Department of An\'alisis Matem\'atico, Facultad de Ciencias, University of M\'alaga, Campus de Teatinos s/n, M\'alaga, 29080-Andaluc\'ia (Spain)}

% \author[auth1]{A. M. Ferreiro-Ferreiro}\ead{ana.fferreiro@udc.es}
% \author[auth1]{J. A. Garc\'ia-Rodr\'iguez}\ead{jose.garcia.rodriguez@udc.es}

\institute{J. G. L\'opez-Salas, M. Su\'arez-Taboada, A. M. Ferreiro-Ferreiro, J. A. Garc\'ia-Rodr\'iguez \at Department of Mathematics, Faculty of Informatics and CITIC, Campus Elvi\~na s/n, 15071-A Coru\~na (Spain), \email{jose.lsalas@udc.es, maria.suarez3@udc.es, ana.fferreiro@udc.es, jose.garcia.rodriguez@udc.es}
\and M. J. Castro \at Department of An\'alisis Matem\'atico, Facultad de Ciencias, University of M\'alaga, Campus de Teatinos s/n, M\'alaga, 29080-Andaluc\'ia (Spain)  \email{mjcastro@uma.es}}

\maketitle
\vspace*{-2cm}
\abstract*{  We present  a novel and general methodology for building second order finite volume  implicit-explicit Runge-Kutta numerical schemes for solving two dimensional financial parabolic PDEs with mixed derivatives.
The methods achieve second order convergence even in the presence of non-regular initial conditions. The IMEX time integrator allows to overcome the tiny time-step induced by the diffusive term in the explicit schemes, also providing accurate and non-oscillatory approximations of the Greeks. 
}
\abstract{ We present  a novel and general methodology for building second order finite volume  implicit-explicit Runge-Kutta numerical schemes for solving two dimensional financial parabolic PDEs with mixed derivatives. The methods achieve second order convergence even in the presence of non-regular initial conditions. The IMEX time integrator allows to overcome the tiny time-step induced by the diffusive term in the explicit schemes, also providing accurate and non-oscillatory approximations of the Greeks. 
}

\vspace{-0.5cm}
\section{Introduction}\label{sec1}
% The goal of this article is to develop second order finite volume numerical schemes for solving option pricing problems, modelled by two dimensional advection diffusion reaction scalar partial differential equations (PDEs). The accuracy of the developed schemes will be measured against alternative and precise spectral Fourier methods. 

The goal of this article is to design  general  second order  numerical schemes for solving financial parabolic PDEs with mixed  derivatives. 
The technique is based in the combination of finite volume (FV) methods, and implicit-explicit (IMEX) Runge-Kutta (RK) time integrators.
To our knowledge, this is the first time that these numerical schemes have been successfully applied to the numerical solution of parabolic financial PDEs.  

% This is a really powerful combination for solving problems in mathematical finance.
% On the one hand, finite volume methods are taylored to treat convective terms, even in highly nonliear hyperbolic PDEs. 
% They are very well stablished, allowing to build schemes with arbitrary order of convergence, even in the presence of non-regular initial or boundary conditions, which is a well-known difficulty in the financial literature.
% On the other hand,  IMEX time integrators allow to overcome the severere time restriction imposed by the spatial semi-discretization of the diffusive terms.  This is also of great importance, because parabolic problems may become intractable from the computational point of view, if explicit methods are applied.  Furthermore, the presented numerical schemes yield excellent approximations of the sesitivities (Greeks), free of the spurious oscillations, which typically appear when low order and no strongly stable time integrators are considered. 

% Last but not least, in order to assess the robustness, accuracy and order of convergence of the proposed finite volume IMEX numerical schemes, the multidimensional Fourier cosine method is implemented. A novel way to truncate the series expansion for basket options is presented. Additionally, an efficient parallel implementation in a multi-GPU environment allows to recover extraordinarily accurate option prices. 

The main  mathematical challenge in the numerical solution of  financial advection-diffusion PDEs arises when the advection term becomes large when compared to the diffusion one. Under this situation, instabilities show up, because the problem becomes more hyperbolic. In order to face these problems, several techniques have been introduced. One way to overcome these instability phenomena is to avoid centered schemes and to consider upwind discretizations of advection terms, thus taking information upstream and not downstream. 
% We refer the reader to \cite{Pironneau05} for further details in the mathematical finance field. 
One technique to properly solve convection-dominated diffusion problems is to consider finite volume discretization methods. 
They are very well stablished, allowing to build schemes with arbitrary order of convergence, even in the presence of large convective terms and non-regular initial or boundary conditions, which are well-documented difficulties in the financial literature.
 Seminal works on applying finite volume method to option pricing problems are due to  Forsyth and Zvan, see  \cite{Zvan-Forsyth-97, ForsythZvan01}, where vertex type finite volume methods were applied for problems in two spatial dimensions. More recently, in \cite{Tadmor2ndOrder2019} the authors propose a Kurganov-Tadmor scheme \cite{KT2000,NT1990} for solving option pricing problems, written in conservative form, with appropriate time integration methods and slope limiters. Finally in \cite{Tadmor2019} the authors apply a third order Kurganov-Levy scheme presented in \cite{KL2000} based in the CWENO reconstructions presented in \cite{LPR1999}. All these works only deal with Black-Scholes PDEs in dimension one or Asian PDE problems. Therefore, they do not solve general two-dimensional pricing problems with mixed derivatives.

In all the previous works the authors use explicit-schemes in time  to solve the stiff systems of ordinary differential equations (ODEs) obtained after performing the spatial semidiscretization. The main disadvantage of these approaches is that tiny time-steps must be used in order to maintain stability. Typically, Von Neumann stability analysis demands $\Delta t \leq {(\Delta x)^2}/{(2 \eta)}$, being $\eta$ the diffusion velocity.
One technique to overcome this restriction is to use implicit-explicit time integrators, see \cite{ASCHER1997, Boscarino2013, Calvo2001,Kennedy2003,Russo05}.
% Pareschi and Russo proposed in \cite{Russo05} implicit-explicit (IMEX) Runge-Kutta method for general hyperbolic systems of conservation laws with stiff relaxation terms. 
The idea is to apply an implicit discretization to the stiff term related to the relaxed source terms, and an explicit one to the nonstiff term. This idea can also be extended to convection diffusion equations, see \cite{russo-finance}.
%, the IMEX time marching scheme was used in mathematical finance only combined with finite differences and finite elements, respectively.
Following these works, our goal is to first semi-discretize in space 2D Black-Scholes equations (having mixed derivatives) by means of second-order finite volume methods, thus properly treating convection terms and non-smooth payoffs. Later, we propose to integrate in time the resulting system of stiff ODEs by means of the second-order IMEX Runge-Kutta time marching scheme. In this way, we will apply an implicit discretization to the diffusion (stiff) part and an explicit one to the convection and source terms (non stiff). 
As a result, $\Delta t$ will only depend on the stability of the convection.
% , being of ${\cal O}(\Delta x)$, instead of ${\cal O}(\Delta x^2)$.
% , i.e, $\Delta t \leq \frac{\Delta x}{\alpha}$, being $\alpha$ the convection velocity. The here proposed numerical methods will provide also good approximations of the Greeks (without oscillations).

% The organization of this paper is as follows. In Section \ref{sec:PDE_Models} we overview a couple of two-dimensional Black-Scholes PDE problems with huge interest in the financial industry.
% Section \ref{sec:numerical-schemes} presents the proposed finite volume IMEX Runge-Kutta scheme for general two-dimensional convection diffusion PDEs written in conservative form. Section \ref{2dcos} is devoted to an alternative Fourier numerical method designed to solve the PDE problems presented in Section \ref{sec:PDE_Models}. In order to obtain a high-accurate reference solution to validate the here proposed finite volume methods, a novel and very efficient multi-GPU implementation of the Fourier method is carried out. In Section \ref{sec:numerical-results} numerical experiments are performed to price several options under the Black-Scholes and Heston models. A thorough inspection of the accuracy and order of convergence of the developed numerical methods is detailed. Finally, conclusions are drawn in Section \ref{conclusions}.

% \subsection{Options on a basket of two assets\label{sec:Basket2D}}
\vspace{-0.5cm}
\section{Option pricing models}

In this work we will focus in two important PDE models in finance: baskets of two assets and the Heston stochastic volatility  model. Both models can be formulated in terms of two dimensional parabolic PDEs.

We consider a basket of two assets, with prices given by $s_1$ and $s_2$. Under the Black-Scholes model, these prices follow the system of stochastic differential equations:
%\begin{align*}
%  \frac{\dx{s_{1t}}}{s_{1t}} &=   (r-d_1) \dx{t} + \dfrac{\sigma_1}{\sqrt{1+\rho^2} } (\dx{W_{1t}} + \rho \dx{W_{2t}}),\\
%  \frac{\dx{s_{2t}}}{s_{2t}} &=   (r-d_2) \dx{t} + \dfrac{\sigma_2}{\sqrt{1+\rho^2} } (\rho \dx{W_{1t}} + \dx{W_{2t}}),  
%\end{align*}
\begin{equation} \label{eq:sdesBS2d}
  % \scalebox{0.5}{ }
  \frac{\dx{s_{it}}}{s_{it}} =   (r-q_i) \dx{t} + \sigma_i \dx{W_{it}}, \quad s_{i0} \mbox{  known},\, i=1,2,
%   &\frac{\dx{s_{2t}}}{s_{2t}} =   (r-q_2) \dx{t} + \sigma_2\dx{W_{2t}}, \quad s_{20} \mbox{ known},
\end{equation} 
where $W$ is a two dimensional correlated Brownian motion. We set 
$W_{1t} = \bar{W}_{1t}$,   $W_{2t} = \rho \bar{W}_{1t} + \sqrt{1-\rho^2} \bar{W}_{2t},$
where $(\bar{W}_{1},\bar{W}_{2})$ is a two dimensional standard Brownian motion and $\rho \in (-1,1)$ is the constant correlation parameter. Besides, $\sigma_1, \sigma_2 \in \R_+$ are the market volatilities of the assets $s_1, s_2$, respectively. Finally, $r$ is the interest rate, and $q_1, q_2$ are the dividend yields of the assets.
The price of an option with maturity $T$ and payoff function $u_T(s_1,s_2)$ is $u(s_{1t},s_{2t},t) = \e{-r(T-t)}\Esp{u_T(s_{1T},s_{2T})},$ where $Q$ is the selected martingale measure.
By applying the two-dimensional Itô formula, we can obtain a backward PDE for the price of the option $u(s_1,s_2,t)$.
% \begin{align}\label{eq:PdeBasket2DBackward}
%     & \dfrac{\partial u}{\partial t}  + \frac12\sigma_1^2s_1^2\frac{\partial^2 u}{\partial s_1^2}+\frac12\sigma_2^2s_2^2\frac{\partial^2 u}{\partial s_2^2} + \rho\sigma_1\sigma_2s_1s_2\frac{\partial^2 u}{\partial s_1\partial s_2} \nonumber\\
%     &\qquad +(r-q_1)s_1 \dfrac{\partial u}{\partial s_1}+(r-q_2)s_2\dfrac{\partial u}{\partial s_2}-ru=0,\qquad t\in[0,T), \nonumber\\
%     & u(s_1,s_2,T) = u_T(s_1,s_2), \qquad s_1,s_2>0.
% \end{align}
Doing a time reversal $\tau = T-t$ change of variable, the following forward PDE is obtained:
{\small
\begin{align}\label{eq:PdeBasket2DForward}
    & \dfrac{\partial u}{\partial \tau }  - \frac12\sigma_1^2s_1^2\frac{\partial^2 u}{\partial s_1^2}-\frac12\sigma_2^2s_2^2\frac{\partial^2 u}{\partial s_2^2} - \rho\sigma_1\sigma_2s_1s_2\frac{\partial^2 u}{\partial s_1\partial s_2} \nonumber\\
    &\qquad -(r-q_1)s_1 \dfrac{\partial u}{\partial s_1}-(r-q_2)s_2\dfrac{\partial u}{\partial s_2}+ru=0,\qquad \tau\in(0,T], \nonumber \\
    & u(s_1,s_2,0) = u_0(s_1,s_2), \qquad s_1,s_2>0.
\end{align}
}
In the sequel, the $\tau$ notation is dropped for simplicity and the forward time is again written as $t$.
In this work we consider the payoff of the arithmetic basket call option,
$
u_0(s_1,s_2)=\max\bigl(\frac{1}{2}(s_1+s_2)-K,0\bigr), 
$
where $K$ is the fixed strike price.

In the Heston model \cite{Heston93}, there is an underlying asset $s$ whose volatility is a stochastic process driven by a second Brownian motion:
{\small
\begin{align}\label{SDE-System-Heston}
&\dx{s_t}=(r-q)s_t \dx{t} +\sqrt{v_t}s_t \dx{W_{1t}},\quad s_0 \mbox{ known,}\nonumber\\
&\dx{v_t}=\kappa(\theta-v_t)\dx{t} +\sigma\sqrt{v_t} \dx{W_{2t}}, \quad v_0 \mbox{ known},
\end{align}
}
where $r, q, \kappa, \theta, \sigma, s_0, v_0$ are constant parameters in $\R_+$ and $W$ is again a two dimensional correlated Brownian motion. Moreover, $r$ is the fixed interest rate, $q$ is the constant dividend yield, $\kappa$ is the mean-reversion speed for the variance, $\theta$ is the mean reversion level for the variance, $\sigma$ is the volatility of the variance (the so-called volatility of the volatility) and $v_0$ is  the initial level of the variance. The process $v$ represents the variance of $s$ and its stochastic differential equation is a version of the square root process described by Cox, Ingersoll and Ross.
The price of a derivative with payoff function $u_T$ is the solution of a backward PDE.
% \begin{align}
% \label{eq:HestonPDEBackward}
% & \dfrac{\partial u}{\partial t}+
% \dfrac{1}{2}s^2v\dfrac{\partial^2 u}{\partial s^2}+
% \rho\sigma s v\dfrac{\partial^2 u}{\partial v\partial s}+
% \dfrac{1}{2}\sigma^2 v \dfrac{\partial^2 u}{\partial v^2}\nonumber \\
% & \qquad +(r-q)s \dfrac{\partial u}{\partial s}+
% \kappa(\theta-v) \dfrac{\partial u}{\partial v}-
% ru=0, \qquad t \in [0,T), \nonumber \\
% &u(s,v,T) = u_T(s), \qquad s>0.
% \end{align}
Once more, performing a time reversal $\tau = T-t$ change of variable, the forward PDE:
{\small
\begin{align}
\label{eq:HestonPDEForward}
& \dfrac{\partial u}{\partial \tau}-
\dfrac{1}{2}s^2v\dfrac{\partial^2 u}{\partial s^2}-
\rho\sigma s v\dfrac{\partial^2 u}{\partial v\partial s}-
\dfrac{1}{2}\sigma^2 v \dfrac{\partial^2 u}{\partial v^2} \nonumber \\
& \qquad -(r-q)s \dfrac{\partial u}{\partial s}-
\kappa(\theta-v) \dfrac{\partial u}{\partial v}+
ru=0, \qquad \tau \in (0,T], \nonumber \\
&u(s,v,0) = u_0(s), \qquad s>0,
\end{align}
}
is obtained. As before, the forward time $\tau$ will be written as $t$.
The payoff function of an European call option is
$
u_0(s) = \max(s-K,0),
$
where $K$ is the strike.

\vspace{-0.5cm}
\section{Finite volume IMEX Runge-Kutta numerical method \label{sec:numerical-schemes}}

Equation \eqref{eq:PdeBasket2DForward} and \eqref{eq:HestonPDEForward}  can be written in the compact conservative form:
% \begin{equation} \label{SistCons}
%     \dfrac{\partial}{\partial t}u(x,y,t)  + \dfrac{\partial f_1}{\partial x} (u)+\dfrac{\partial f_2}{\partial y}  (u)=\dfrac{\partial g_1}{\partial x} (u_x,u_y)+\dfrac{\partial g_2}{\partial y} (u_x,u_y)+h(u),
% \end{equation}
% or
{\small
\begin{equation} \label{SistConsCompacto}
    \dfrac{\partial}{\partial t}u(x,y,t)  + \div{\mbF}(u)=\div{\mbG}(\nabla u)+h(u),
\end{equation}
}
with $\mbF(u)=(f_1(u),f_2(u))$ and $\mbG(\nabla u)=(g_1(\nabla u),g_2(\nabla u))$. 
% More precisely, $\mbF:\R^3\rightarrow\R^2$,  $\mbG:\R^4\rightarrow\R^2$ and $h:\R^3\rightarrow\R$ are functions of $x,y,u,\nabla u$, although in order to keep the notation simple some of these dependencies were omitted.
The numerical solution of equation \eqref{SistConsCompacto} using a finite volume scheme  is difficult due to the diffusive terms.  The obtained semidiscrete scheme is a stiff ODE system:
{\small
\begin{equation}\label{eq:stiffedo}
\dfrac{\partial U}{\partial t}  +  E (U)=I (U),
\end{equation}
}
where $U=U(t)\in \R^N$ and $E,I:\R^N\to \R^N$, being $E$ the non-stiff term and $I$  the stiff part. 
IMEX  time-marching schemes play a major rule in the treatment of these stiff systems.
In this article we have considered the implicit-explicit (IMEX) Runge-Kutta time discretization numerical scheme proposed in \cite{Russo05}.  

% Both parts will be handled simultaneously with the same IMEX solver. 

% The rest of the section is organized as follows. In Section \ref{sec:spatialDiscretization} we describe the space discretization obtained by finite volume schemes. Section \ref{sec:timeDiscretization} is devoted to the IMEX Runge-Kutta scheme applied to the stiff system of differential equations \eqref{eq:stiffedo} obtained by the finite volume space discretization.

\vspace{-0.5cm}
\subsection{Finite volume space discretization} \label{sec:spatialDiscretization}

Let $\Delta x$, $\Delta y$ be the mesh length in each spatial direction.
We define the grid points $x_i=i\Delta x$, $y_j=j\Delta y$, $i,j \in \mathbb{Z}$. Let $x_{i+1/2} = x_i + \frac12\Delta x$ and $y_{j+1/2} = y_j + \frac12\Delta y$. We consider rectangular finite volumes $V_{ij}= [x_{i -1/2} , x_{i +1/2} ]\times[y_{j -1/2} , y_{j +1/2} ] $, where $(x_i,y_j)$ is the center of the finite volume $V_{ij}$. The area of $V_{ij}$ will be denoted as $\lvert V_{ij}\rvert$, i.e.  $\lvert V_{ij}\rvert = \Delta x\Delta y$. Besides, $\Gamma_{ij}$ represents the boundary of $V_{ij}$. Let $\bar{u}_{ij}=\frac{1}{\lvert V_{ij}\rvert}\int_{V_{ij}} u(x,y,t) \dx{x}\dx{y}$ be the volume averages that will be the unknowns of the problem . 
% The here described finite volume grid is sketched in Figure \ref{fig:Stencil}.
% \begin{figure}
% \centering
% \includegraphics[width=8cm]{./images/grafica_2D_v2.pdf}
% \caption{Finite volume stencil.}\label{fig:Stencil}
% \end{figure}
Integrating equation \eqref{SistConsCompacto} in space on $V_{ij} $ and dividing by $\lvert V_{ij} \rvert$ we obtain the semi-discrete equation
{\small
\begin{align}
\label{eq:integrada_vols}
\dfrac{\dx{\bar{u}_{ij}}}{\dx{t}}
 =-\dfrac{1}{ \lvert V_{ij}\rvert} \int_{V_{ij}}  \div{\mbF}(u) \,\dx{x}\dx{y} +\dfrac{1}{\lvert V_{ij} \rvert} \int_{V_{ij}}  \div{\mbG}(u_x,u_y) \,\dx{x}\dx{y}
%  + \overline{r(u)}_i .
+\dfrac{1}{ \lvert V_{ij} \rvert} \int_{V_{ij}} h(u) \, \dx{x}\dx{y}.
\end{align}
}
Applying the divergence theorem:
%  we can rewrite the first two volume integrals of \eqref{eq:integrada_vols} as line integrals
{\small
\begin{align}
\label{eq:integragada_lineIntegrals}
\dfrac{\dx{ \bar{u}_{ij}}}{\dx{t}}
 =-\dfrac{1}{ \lvert V_{ij} \rvert} \oint_{\Gamma_{ij}}  {\mbF}(u)\cdot{\mb{ n} } \,d\gamma 
+\dfrac{1}{ \lvert V_{ij} \rvert}  \oint_{\Gamma_{ij}}  {\mbG}(u_x,u_y)\cdot{\mb{ n} } \,d\gamma 
+\dfrac{1}{ \lvert V_{ij} \rvert} \int_{V_{ij}} h(u) \, \dx{x}\dx{y}.
\end{align}
}
Therefore, in order to convert \eqref{eq:integragada_lineIntegrals} into a numerical scheme, we have to
approximate on the right hand side of this equation with functions of $\{\bar{u}_{ij}(t)\}$. 
% In the following two content blocks, the numerical treatment of the source and advection terms (Section \ref{sec:numSemiDiscExplicit}) and the diffusion ones (Section \ref{sec:numSemiDiscImplicit}) is discussed. 
The source and convective terms will be treated explicitly in time, while the diffusion part will be managed implicitly.
% \subsubsection{Explicit part: advection-reaction space discretization} \label{sec:numSemiDiscExplicit}
For the advective part we get
{\footnotesize
\begin{align}
% \label{eq:semidiscretoespacio}
 \oint_{\Gamma_{ij}} {\mbF}(u){\mb{ n} } =& \int_{\Gamma_{i+1/2},j} f_1(u) \dx{\gamma_j} + \int_{\Gamma_{i-1/2,j}} f_1(u) \dx{\gamma_j}+ \int_{\Gamma_{i,j+1/2}} f_2(u)  \dx{\gamma_i} + \int_{\Gamma_{i,j-1/2}} f_2(u )  \dx{\gamma_i}, \label{eq:explicitPartAdvection}
\end{align} 
}
where for second order schemes each line integral can be approximated by means of  midpoint quadrature rule as 
$
   \int_{\Gamma_{i\pm 1/2,j}} f_1(u) \dx{\gamma_j}\approx \pm \Delta y f_1(u(x_{i\pm 1/2},y_j)),
%    \int_{\Gamma_{i,j\pm1/2}} f_2(u)  \dx{\gamma_i} &\approx  \pm \Delta x f_2(u(x_i,y_{j\pm1/2})).
$
and analogously for $\int_{\Gamma_{i,j\pm1/2}} f_2(u)  \dx{\gamma_i}$.
At this point the unknown function $u$ is reconstructed by a piecewise polynomial using the volume averages $\{\bar{u}_{ij} (t)\}$. More precisely, starting from $\{\bar{u}_{ij}\}$, we compute a piecewise polynomial reconstruction
$
\mathcal{R}(x,y) = \sum_{i,j}
 P_{ij} (x,y) \mathbf{1}_{ij} (x,y),
$
where $P_{ij}$ is a suitable polynomial satisfying some accuracy and non oscillatory property, and $\mathbf{1}_{ij}$ is the indicator function of the volume $V_{ij}$. 
Second order numerical schemes can be obtained by means of piecewise linear polynomials, although higher order schemes can be obtained by polynomials of higher order. Here we consider the natural extension of MUSCL reconstruction to 2D Cartesian grids (see \cite{BramVanLeer}).
The flux functions at the midpoints of the boundaries of the volumes can be computed by using a
suitable numerical flux function, consistent with the analytical flux.
Hereafter we detail the approximation of $f_1 (u(x_{i +1/2},y_j))$, being the other fluxes approximated in the same way:
$
f_1 (u(x_{i +1/2},y_j )) \approx \mathcal{F}_1 (u_{i+1/2,j}^- , u_{i +1/2,j}^+ ).
$
The values $u^\pm_{i +1/2,j}$ are obtained from the reconstruction as
$
u^\pm_{i +1/2,j}
 =
 \lim_{x\to x^\pm_{i +
 1/2}}
 \mathcal{R}(x,j\Delta y),
 $
 with $x$ in a normal line to the boundary $\Gamma_{i+1/2,j}$. For example the left reconstructed value at the edge $\Gamma_{i+1/2,j}$ is
$
u^-_{i + 1/2,j}= \bar{u}_{i,j}+\Delta x/2\, u'_{i,j},
$
where the slope $u'_{i,j}$ is a first order approximation of the space derivative in the $X$ direction of $u(x,y,t)$ at point $(x_i, y_j)$ at every time $t$.
This slope must satisfy the TVD property and thus we must use slope limiters. In our case we use the minmod limiter.
% , where the slope is given by
% $$
% u'_{ij}=\dfrac{1}{\Delta x} {\rm minmod} ( \bar{u}_{i,j}-\bar{u}_{i-1,j},\bar{u}_{i+1,j}-\bar{u}_{i,j}),
% $$
% $$
% {\rm minmod}(a,b)=
% \begin{cases}
% \min (a,b) & \text{if } a,b>0,\\
% \max(a,b) & \text{if } a,b<0,\\
% 0 & \text{otherwise\,. } 
% \end{cases}
% $$
All the calculations performed in this article have been carried out using the Local Lax Friedrich numerical flux.
% $$
% \mathcal{F}_1(u^-_{i +1/2,j},u^+_{i +1/2,j})=\dfrac{1}{2}(f_1(u^-_{i +1/2,j})+f_1(u^+_{i +1/2,j}))- \dfrac{1}{2}\alpha(u^+_{i +1/2,j}-u^-_{i +1/2,j}),
% $$
% where $\alpha=\left\lvert f'_1\left((u^-_{i +1/2,j}+u^+_{i +1/2,j})/{2}\right)\right\rvert$.
% Therefore, the line integral \eqref{eq:explicitPartAdvection} of the convective term, is finally approximated as
% \begin{align}
% \oint_{\Gamma_{ij}}  {\mbF}(u)\cdot{\mb{ n} } \,d\gamma \approx& \Delta y \left( \mathcal{F}_1 (u^-_{i +1/2,j} , u^+_{i+1/2,j} ) -\mathcal{F}_1 (u^-_{ i -1/2,j} , u^+_{i-1/2,j} ) \right)+\notag \\  
%  & \Delta x \left ( \mathcal{F}_2 (u^-_{i, j +1/2},u^+_{i,j+1/2} ) -\mathcal{F}_2 (u^-_{i,j -1/2},u^+_{i,j-1/2} ) \right),
% \label{eq:semidiscretospaceConvection} 
% \end{align}
% being $\mathcal{F}_1$ and $\mathcal{F}_2$ the numerical fluxes of the physical flux functions $f_1$ and $f_2$, respectively.
Finally, the volume integral of the source term is discretized using the midpoint  quadrature rule.
% \begin{align}\label{eq:semidiscretospaceReaction}
% \int_{V_{ij}}  h(u) \dx{x}\dx{y}\approx \lvert V_{ij}\rvert h(\bar{u}_{ij}).
% \end{align}

% \subsubsection{Implicit part: diffusion space discretization} \label{sec:numSemiDiscImplicit}

The diffusion part is discretized implicitly:
{\small
\begin{align}
 \oint_{\Gamma_{ij}} {\mbG}(u_x,u_y)\cdot{\mb{ n} } \,d\gamma=& \int_{\Gamma_{i+1/2},j} g_1(u_x,u_y) \dx{\gamma_j} + \int_{\Gamma_{i-1/2,j}} g_1(u_x,u_y) \dx{\gamma_j}  + \notag\\
 & \int_{\Gamma_{i,j+1/2}} g_2(u_x,u_y)  \dx{\gamma_i} + \int_{\Gamma_{i,j-1/2}} g_2(u_x,u_y)  \dx{\gamma_i}. \label{eq:implicitPartDiffusion}
\end{align} 
}
Before approximating each one of these line integrals, a suitable approximation of the partial derivatives $u_x$ and $u_y$ has to be built for the volume $V_{ij}$. With this aim, we build the second order Lagrange interpolating polynomial of $u$ centered in the volume $V_{ij}$. Let $L_{ij}^u$ be that polynomial. Considering the  nine nodes of the dual mesh $\{x_{i+k},y_{j+l}\},\, k,l=-1,0,1$,  and the averaged values $\{\bar{u}_{i+k,j+l}\}$ of the solution at each volume ${V_{i+kj+l}}$, this polynomial is given by
{\small
$
L_{ij}^{\bar{u}}(x,y)=\sum_{k,l=-1}^1 \bar{u}_{i+k,j+l} \, \ell_{i+k}(x)\ell_{j+l}(y), 
$
}
where $\ell_{i+k}$ and $\ell_{j+l}$ are the one dimensional Lagrange polynomial basis, i.e:
{\small
$$
 \ell_{i+k}(x) = \prod_{\substack{p=i-1 \\ p\neq i+k}}^{i+1} \frac{x-x_p}{x_{i+k}-x_p}, \qquad \ell_{j+l}(y) = \prod_{\substack{q=j-1 \\ q\neq j+l}}^{j+1} \frac{y-y_q}{y_{j+l}-y_q}.
$$
}
Therefore, we use the approximations  $u_x\approx \partial_x L_{ij}^{\bar{u}}$ and $u_y\approx \partial_y  L_{ij}^{\bar{u}}$. 
% The computation of these approximations of the partial derivatives of the solution is straightforward, since it just involves the computation of derivatives of one dimensional Lagrange polynomial basis.
For example, the line integrals in equation \eqref{eq:implicitPartDiffusion} are approximated, in the $X$ direction as
{\small
\vspace{-0.25cm}
\begin{align*}
%    \int_{\Gamma_{i\pm 1/2,j}}  & g_1( u_x, u_y)  \dx{\gamma_j}  \approx  
   \int_{\Gamma_{i\pm 1/2,j}}  g_1\Bigl( \partial_x L_{ij}^{\bar{u}}, \partial_y L_{ij}^{\bar{u}}\Bigr)  \dx{\gamma_j}\approx
& \pm \Delta y\,  g_1\Bigl( \partial_x L_{ij}^{\bar{u}}(x_{i\pm 1/2},y_j  ), \partial_y L_{ij}^{\bar{u}}(x_{i\pm 1/2},y_j  )\Bigr).
% \\
% \\
%    \int_{\Gamma_{i,j\pm 1/2}}  & g_2( u_x, u_y)  \dx{\gamma_i}  \approx  \int_{\Gamma_{i,j\pm 1/2}}  g_2\Bigl( \partial_x L_{ij}^{\bar{u}}, \partial_y L_{ij}^{\bar{u}}\Bigr)  \dx{\gamma_i}\approx
% \\
% & \pm \Delta x\,  g_2\Bigl( \partial_x L_{ij}^{\bar{u}}(x_i, y_{j\pm 1/2}  ), \partial_y L_{ij}^{\bar{u}}(x_i, y_{j\pm 1/2} )\Bigr).
\end{align*}
}
% and analogously for $ \int_{\Gamma_{i,j\pm 1/2}}  g_2( u_x, u_y)  \dx{\gamma_i}$.

\vspace{-1cm}
\subsection{IMEX Runge-Kutta time discretization} \label{sec:timeDiscretization}

An IMEX Runge-Kutta scheme consists of applying an implicit discretization to the diffusion terms (stiff terms) and an explicit one to the convective and source terms (non stiff terms), see \cite{ASCHER1997,Boscarino2013, Calvo2001,Kennedy2003,Russo05}. When applied to system \eqref{eq:stiffedo} it takes the form
{\small
\begin{align}
U^{(k)} &= U^n -\Delta t \sum_{l=1}^{k-1} \tilde{a}_{kl} {E}(U^{(l)})+\Delta t \sum_{l=1}^{s} a_{kl} {I}(U^{(l)}), \label{eq:RKIMEX_steps}\\
U^{n+1} &= U^n- \Delta t \sum_{k=1}^{s} \tilde{\omega}_{k} E(U^{(k)})+\Delta t \sum_{k=1}^{s} \omega_{k}I(U^{(k)}), \label{eq:RKIMEX_final}
\end{align}
}
where $U^n=(\bar{u}_{ij}^n)$, $U^{n+1}=(\bar{u}_{ij}^{n+1})$ are the vectors of unknowns volume averages at  times $t^n$ and $t^{n+1}$, and $U^{(k)}$ and $U^{(l)}$ are the vector of unknowns at the stages $k,l$ of the IMEX Runge-Kutta scheme. The matrices $\bar{A} = (\tilde{a}_{ij})$, $\tilde{a}_{kl}=0$ for $l\geq k$ and $A = (a_{kl})$ are $s\times s$ matrices such that the scheme is explicit in $E$ and implicit in $I$. The coefficient vectors $\tilde{w} = (\tilde{w}_1,\ldots,\tilde{w}_s)$ and $w = (w_1,\ldots,w_s)$ complete the IMEX Runge-Kutta scheme.
One should consider diagonally implicit Runge-Kutta (DIRK) schemes for the implicit terms ($a_{kl} = 0$, for $l > k$) in order to solve efficiently the algebraic equations corresponding to the implicit part of the discretization at each time step. 
% Besides, the use of a DIRK scheme for the implicit part $I$ guarantees that $F$ is always evaluated explicitly.
IMEX Runge-Kutta schemes can be represented by a double \textit{tableau} in the usual Butcher
notation,

\begin{center}
\begin{tabular}{c|c}
$\tilde{c}$  & $\tilde{A}$ \\
\hline
&  $\tilde{\omega}$ 
\end{tabular}
\quad
\begin{tabular}{c|c}
$c$  & $A$ \\
\hline
         &  $\omega$ 
\end{tabular},
\end{center}
where the coefficients vectors $\tilde{c}=(\tilde{c}_1,\ldots,\tilde{c}_s)^\mathsf{T}$ and $c=(c_1,\ldots,c_s)^\mathsf{T}$ used for the treatment of non autonomous systems, are given by the
usual relation
$ \tilde{c}_k=  \sum_{l=1}^{k-1} \tilde{a}_{kl}$, \mbox{$c_k=\sum_{l=1}^{k} {a}_{kl}.$}
The PDEs \eqref{eq:PdeBasket2DForward} and \eqref{eq:HestonPDEForward} we will be dealing with in this article, after the space discretization, give rise to a system of non-autonomous ODEs. Therefore, the constants $\tilde{c}_k$ and $c_k$ were omitted in the time marching scheme \eqref{eq:RKIMEX_steps}-\eqref{eq:RKIMEX_final}.
In this work we will consider the second order IMEX-SSP2(2,2,2) (L-stable scheme, see \cite{Russo05}), whose \textit{tableaus} for the explicit (left) and implicit (right) parts are
{\small
\begin{center}
\begin{tabular}{c|cc}
0  & 0 & 0 \\
1  & 1 & 0 \\
\hline
& $1/2$ & $1/2$  
\end{tabular}
\quad
\begin{tabular}{c|cc}
$\gamma$  & $\gamma$ & 0 \\
$1-\gamma$  & $1-2\gamma$ & $\gamma$ \\
\hline
& $1/2$ & $1/2$  
\end{tabular}
\quad 
$\gamma=1-\dfrac{1}{\sqrt{2}}.$ 
\end{center}
}
\noindent The selection of this second-order IMEX scheme is only to match it to the finite
volume framework and to show that globally the scheme works well in the context of parabolic
PDEs finance. Nevertheless,  many other second-order IMEX schemes with $s = 2$,
and with the same properties (stiff problems, L-stable) could be applied in this context as well, see \cite{ASCHER1997} for example.

\vspace{-1.cm}
\section{Numerical experiments \label{sec:numerical-results}}

In this section, several numerical experiments are developed to assess the accuracy and performance of the new numerical methods proposed in previous sections for the discussed two dimensional problems. 
Experiments are carried out on several option pricing problems. In subsection \ref{sec:NumExpBasket} basket call options over two underlyings following Black-Scholes model are priced. Vanilla call options under Heston stochastic volatility model are priced in subsection \ref{sec:NumExpHeston}. 
Each subsection is organized as follows. 
We start by writing the PDE model in conservative form. Then, we solve the PDE problems using the discussed finite volume IMEX Runge-Kutta numerical method. Both convection-dominated and diffusion-dominated scenarios are considered. 
% As a result, surface prices are presented. In addition, plots for the numerical derivatives of the solution are presented. 
% More precisely, delta and gamma Greeks are shown. 
All pictures, tables and errors in this whole section will be computed considering the numerical solution at the last time step where $t=T$, i.e, we approximate numerically the option prices ``today''.
Later, reference solutions for the studied pricing problems are accurately computed at $t=T$ by means of the COS Fourier method \cite{fang:oosterlee:08}.
%  explained in Section \ref{2dcos}. 
 This allow us to compute $L_1$ and $L_\infty$ error norms, together with the relative and mean absolute errors.
% given by:
% \begin{align*}
%     L_1 \mbox{ error } &= \Delta x \Delta y\sum_{i,j}  \lvert \bar{u}_{i,j}-u(x_i,y_j) \rvert, \\
%     L_\infty \mbox{ error } &= \max_{i,j} \lvert \bar{u}_{i,j}-u(x_i,y_j) \rvert, \\    
%     L_\infty \mbox{ relative error } &= \dfrac{L_\infty \mbox{ error }}{\max_{i,j} \lvert u(x_i,y_j) \rvert},\\
%     \mbox{Mean absolute error} &=\underset{i, j}{ \mbox{mean} } \lvert \bar{u}_{i,j}-u(x_i,y_j) \rvert.
% \end{align*}
% Besides, absolute errors of the finite volume IMEX Runge-Kutta method against the COS Fourier method are plotted in order to quickly comprehend the quality of the approximation and the regions with higher errors. 
On top of that, orders of convergence are validated in-depth. Both explicit and IMEX Runge-Kutta time integrators were implemented over the same explained space discretization. In the explicit case second order Heun's ODE solver was considered. Finite volume with both time marching schemes achieves the announced second order accuracy in all regimes. IMEX and explicit methods yield similar results in terms of accuracy and convergence order. Nevertheless, the IMEX solver can take advantage of much larger times steps, thus offering much better performance in terms of execution times. The explicit method dramatically slows down the computation due to the presence of diffusive terms.
% Finally, each subsection is closed with a table with some reference option prices. Since the here considered option prices are not known is closed form, an interested reader may find this table useful, particularly those developers of alternative option price calculators.
The designed algorithms were implemented using C++ (GNU compiler 9.3.0) and ran in a machine with an AMD Zen3 5950X processor. All codes were compiled using double precision. 
% Besides, for linear algebra operations Eigen library \cite{eigenweb} was considered. 
In this article a CFL of $0.5$ is taken into account in the stability conditions.

\begin{figure}[!htb]
  \centering
  \subfigure[Test 1] {\includegraphics[height=3.5cm]{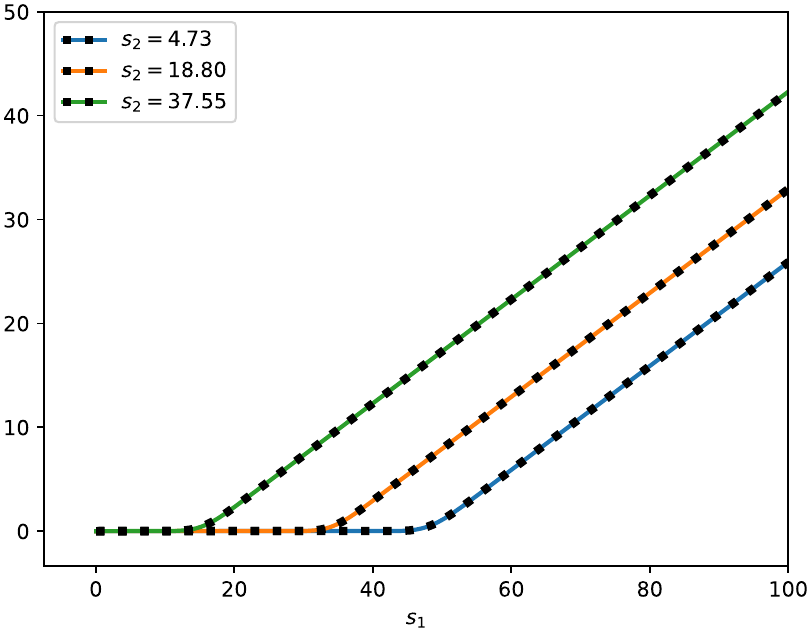}}
  \subfigure[Test 2] {\includegraphics[height=3.5cm]{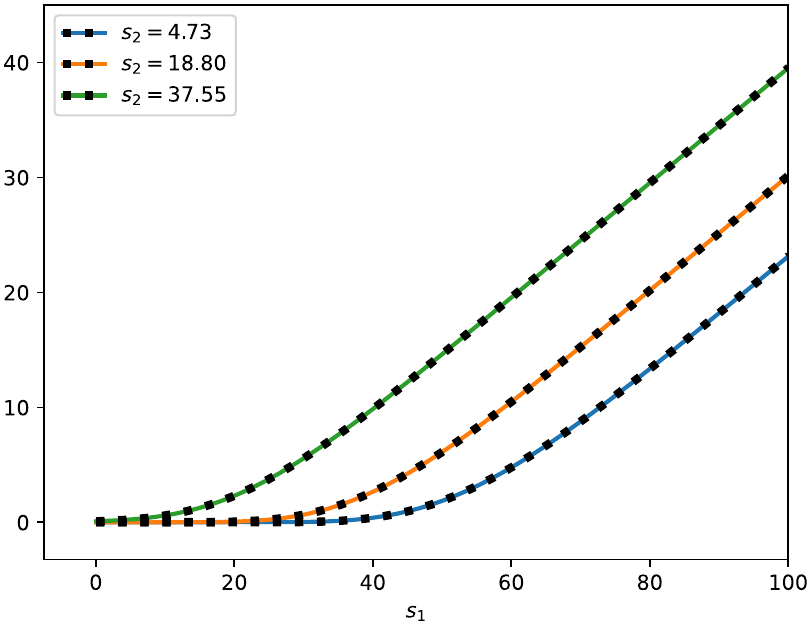}}
  \caption{Cuts of price surfaces. Numerical solution, continuous; COS solution with squares.}
  \label{fig:basket-prices-cuts}
  \end{figure}

  \vspace{-1.5cm}
\subsection{Options on a basket of two assets\label{sec:NumExpBasket}}

In this experiment we solve the basket European option pricing problem considering two underlyings. The PDE model \eqref{eq:PdeBasket2DForward} can be written in the conservative form \eqref{SistConsCompacto}:
{\small
\begin{equation} \label{SistConsBS}
    \dfrac{\partial u}{\partial t}  + \dfrac{\partial f_1}{\partial s_1} (u)+\dfrac{\partial f_2}{\partial s_2}  (u)=\dfrac{\partial g_1}{\partial s_1} (u_{s_1},u_{s_2})+\dfrac{\partial g_2}{\partial s_2} (u_{s_1},u_{s_2})+h(u),
\end{equation}
}
where the functions $f_1, f_2, g_1, g_2$ and $h$ are given by:
\vspace{-0.3cm}
{\small
\begin{align*}
f_i(u)&=(\sigma_i^2-r+q_i)s_iu(s_1,s_2,t)+\frac{\rho}{2}\sigma_1\sigma_2s_i u(s_1,s_2,t), \\
% f_2(u)&=(\sigma_2^2-r+q_2)s_2u(s_1,s_2,t)+\frac{\rho}{2}\sigma_1\sigma_2s_2 u(s_1,s_2,t), \\
g_1(u_{s_1},u_{s_2})&=\frac12\sigma_1^2s_1^2u_{s_1}(s_1,s_2,t)+\frac\rho2\sigma_1\sigma_2s_1s_2 u_{s_2}(s_1,s_2,t),\\
g_2(u_{s_1},u_{s_2})&=\frac12\sigma_2^2s_2^2 u_{s_2}(s_1,s_2,t)+\frac\rho2\sigma_1\sigma_2s_1s_2 u_{s_1}(s_1,s_2,t),\\
h(u)&=( \sigma_2^2+\rho \sigma_1\sigma_2+\sigma_1^2+q_1+q_2-3r ) u(s_1,s_2,t).
\end{align*}
}
% The initial condition $u_0$, being equation \eqref{eq:arithmeticPayoff} for the arithmetic basket call option, was averaged on each volume. 
% Indeed, the initial volume averages are set as
% $${\bar{u}^0}_{i,j} = \dfrac{1}{\lvert V_{ij} \rvert} \int_{V_{ij}} u_0(s_1,s_2)\dx{s_2}\dx{s_1},$$
% where the mid-point quadrature formula is used to approximate the integral.

We perform two basket option pricing tests with the market parameters given by $q_1=0$, $q_2=0$, $\rho=0.5$, $T=0.25$, $K=30$, and $\sigma_1=\sigma_2=0.1$, $r=0.5$ for Test 1, and 
$\sigma_1=\sigma_2=0.5$, $r=0.1$ for Test 2.
The first set of parameters, denoted as Test 1, stands for a convection-dominated regime. 
Although nowadays this setup with an interest rate $r=0.5$ is financially unrealistic, it is useful as a stress-test. The second group of parameters for Test 2, represents a diffusion-dominated setting. In both scenarios, in order to avoid noise coming from the artificial boundary conditions, the computational space domain was considered as $(s_1,s_2)\in[0,5K]\times[0,5K]$. 

% \begin{table}[h]
% \begin{center}
% \caption{Market data for basket options under Black-Scholes model.}\label{tb-Basket}%
% \begin{tabular}{c|c|c|c|c|c|c|c|c}
%   &$\sigma_1$ & $\sigma_2$ & $r$  & $q_1$ & $q_2$  & $\rho$ & $T$& $K$\\
% \hline
% Test 1& $0.1$ & $0.1$  & $0.5$ & \multirow{2}{*}{$0.00$} & \multirow{2}{*}{$0.00$} & \multirow{2}{*}{$0.5$} & \multirow{2}{*}{$0.25$} & \multirow{2}{*}{$30$}\\

% Test 2 &$0.5$ & $0.5$  & $0.1$ &  &  &  &  & \\
% \end{tabular}
% \end{center}
% \end{table}

% Surfaces of basket option prices at $t=T$ for Test 1 and Test 2 are shown in Figure \ref{fig:basket-prices}. 

Cuts of the surfaces, with planes parallel to $s_2=0$, computed with the finite volume IMEX 
Runge-Kutta solver, for basket option prices at $t=T$ are shown in Figure \ref{fig:basket-prices-cuts}, 
along with COS solutions. 
% All these plots were drawn using the solution obtained with the finite volume IMEX 
% Runge-Kutta solver. 
% Moreover, first numerical derivative (delta) and second numerical derivative (gamma) 
% of the numerical approximation of the option prices are shown in Figures \ref{fig:basket-delta} and 
% \ref{fig:basket-gamma}, respectively. 
The finite volume IMEX Runge-Kutta numerical scheme offers high 
resolution approximation of basket option values. 
% Furthermore, it provides accurate and non-oscillatory 
% approximations of the Greeks.
% \begin{figure}[!htb]
% \centering
% \subfigure[Test 1] {\includegraphics[height=4cm]{./Test_Basket/precio_basket_convectivo}}
% \subfigure[Test 2] {\includegraphics[height=4cm]{./Test_Basket/precio_basket_difusivo}}
% \caption{Basket option prices for Test 1 (left) and Test 2 (right) at $t=T$.}
% \label{fig:basket-prices}
% \end{figure}

% \begin{figure}[!htb]
% \centering
% \subfigure[Test 1] {\includegraphics[height=4cm]{./Test_Basket/delta_basket_convectivo}}
% \subfigure[Test 2] {\includegraphics[height=4cm]{./Test_Basket/delta_basket_difusivo}}
% \caption{Basket option deltas for Test 1 (left) and Test 2 (right) at $t=T$.}
% \label{fig:basket-delta}
% \end{figure}
% \begin{figure}[!htb]
% \centering
% \subfigure[Test 1] {\includegraphics[height=4cm]{./Test_Basket/gamma_basket_convectivo}}
% \subfigure[Test 2] {\includegraphics[height=4cm]{./Test_Basket/gamma_basket_difusivo}}
% \caption{Basket option gammas for Test 1 (left) and Test 2 (right) at $t=T$.}
% \label{fig:basket-gamma}
% \end{figure}
At this point we compute the reference option prices by using the described COS method. 
For a mesh of size $N_1\times N_2 = 1600\times1600$, i.e., a mesh with 1600 discretization points in each 
space direction.
% , the errors of the approximation obtained with the finite volume IMEX Runge-Kutta method 
% are collected in Table \ref{tb-Basket_errors}. 
% Figure \ref{fig:error_basket} presents contour plots for the 
% surface defined by the absolute errors in the option prices at $t=T$.
The here developed finite volumes (explicit or IMEX) 
schemes offer high-resolution approximations even at regions of discontinuities and non-smoothness in the initial condition. 
% Therefore, there is no need to apply smoothing methods like the Rannacher technique.
% \begin{table}[!htb]
% \begin{footnotesize}
% \begin{center}
% \caption{Finite volume IMEX Runge-Kutta numerical errors against reference basket option prices computed with COS method at $t=T$.}
% \begin{tabular}{|c|c|c|}
% \hline
% \multicolumn{3}{|c|}{Test 1}\\
% \hline
%   $L_\infty  \mbox{ error }$ & $L_\infty \mbox{ relative error }$ & $\mbox{Mean absolute error}$ \\
% \hline
% $1.0245\times 10^{-4}$ & $1.3976\times 10^{-6}$ & $3.5079\times 10^{-6}$  \\
% \hline
% \hline
% \multicolumn{3}{|c|}{Test 2}\\
% \hline
%   $L_\infty  \mbox{ error }$ & $L_\infty \mbox{ relative error }$ & $\mbox{Mean absolute error}$ \\
% \hline
% $2.3406\times 10^{-5}$ & $1.9439\times 10^{-7}$ & $2.7099\times 10^{-6}$  \\
% \hline
% \end{tabular} 
% \label{tb-Basket_errors} 
% \end{center}
% \end{footnotesize}
% \end{table}
% \begin{figure}[!htb]
% \centering
% \subfigure[Test 1] {\includegraphics[height=4cm]{./Test_Basket/error_basket_convectiva_r05_sigma01}}
% \subfigure[Test 2] {\includegraphics[height=4cm]{./Test_Basket/error_basket_difusiva_r01_sigma05_dom150_v2}}
% \caption{Contour plots for the surface of absolute errors of basket options prices for Test 1 (left) and Test 2 (right) at $t=T$.}
% \label{fig:error_basket}
% \end{figure}
Table \ref{tb-Basket_diffusive} records $L_1$ errors and orders of convergence for Test 2. The errors and orders of convergence are shown for both the IMEX and explicit finite volume numerical schemes. Both schemes achieve second-order accuracy in the $L_1$ norm. Additionally, these tables present the time steps and the execution times for the two time integrators. As said before, IMEX method is able to converge using much larger time steps than the explicit one. The stability condition of the explicit scheme requires extremely small time steps, thus making the method useless in practice for refined grids in space. Therefore, IMEX offers much better performance in terms of execution times, and allows us to solve the PDE problems with refined meshes in space. In fact, IMEX method is able to run for grids finer than $800\times 800$ grid, while the explicit method is not in reasonable computational times. 
% In Figure \ref{fig:curveEfficiency}, the natural logarithms of $L_1$ errors and execution times of Table \ref{tb-Basket_diffusive} are shown for both the IMEX and explicit numerical schemes. 
IMEX clearly outperforms the explicit method in this scenario with diffusion dominance, which is the typical situation in finance. Although for convective dominated problems in coarse grids both time marching schemes perform similarly, as soon as the mesh is refined IMEX is the only practical choice.

\begin{table}[!h]
% \begin{center}
{\scriptsize
\vspace{-0.5cm}
\begin{tabular}{|c||c|c|c|c|}
\hline
& \multicolumn{4}{c|}{IMEX}  \\
\hline
$N_1\times N_2$  & $L_1$ \text{ error} &  \text{ Order}& $\Delta t$ & Time (s) \\
\hline
$25\times 25$ & $9.6620\times 10^{1}$ & $--$ & $4.71\times 10^{-2}$ & $4.5\times 10^{-3}$ \\
$50\times 50$ &  $2.5178\times 10^{1}$ & $1.94$& $2.35\times 10^{-2}$ & $3.3\times 10^{-2} $ \\
$100\times 100$ &  $6.4828\times 10^{0}$ & $1.95$& $1.18\times 10^{-2}$ & $1.7 \times 10^{-1} $ \\
$200\times 200$ &  $1.6209\times 10^{0}$ & $2.00$ & $5.89\times 10^{-3}$ & $1.2\times 10^{0} $ \\
$400\times 400$ &  $3.9419\times 10^{-1}$ & $2.03$ & $2.94\times 10^{-3}$ & $9.8\times 10^{0}$ \\
$800\times 800$ &  $7.9229\times 10^{-2}$ & $2.31$ & $1.47\times 10^{-3}$ & $8.5\times 10^{1}$ \\
\hline
\end{tabular} 
}
{\scriptsize
\begin{tabular}{|c|c|c|c|}
    \hline
& \multicolumn{3}{c|}{Explicit}  \\
\hline
 $L_1$ \text{ error} &  \text{ Order}& $\Delta t$ & Time (s) \\
\hline 
 $9.0224\times 10^{1}$ & $--$ & $1.42\times 10^{-3}$ & $1.7\times 10^{-2}$\\
 $2.3440\times 10^{1}$ &$1.94$ &  $3.56\times 10^{-4}$ & $1.2\times 10^{-1}$\\
 $5.9498\times 10^{0}$ &$1.97$ & $8.89\times 10^{-5}$ & $1.8\times 10^{0}$\\
 $1.4834\times 10^{0}$ &$2.00$ & $2.22\times 10^{-5}$ & $3.0\times 10^{1}$\\
 $3.5473\times 10^{-1}$ &$2.06$ & $5.56\times 10^{-6}$ & $4.9 \times 10^{2}$\\
 $7.1095\times 10^{-2}$ &$2.31$ & $1.39\times 10^{-6}$ & $ 7.9 \times 10^{3}$ \\
\hline
\end{tabular} 
}
% \end{center}
\caption{$L_1$ errors and orders of convergence of the IMEX and explicit methods for Test 2.}
\label{tb-Basket_diffusive}
\end{table}
\vspace{-1.cm}
\subsection{Heston model\label{sec:NumExpHeston}}
\begin{figure}[!htb]
\vspace{-0.5cm}
\centering
\subfigure[Test 3] {\includegraphics[height=3.5cm]{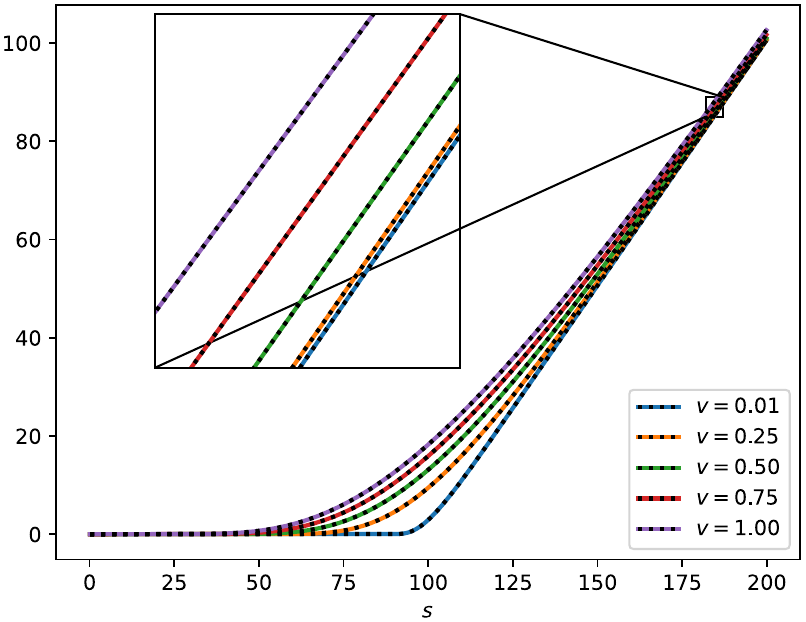}}
\subfigure[Test 4] {\includegraphics[height=3.5cm]{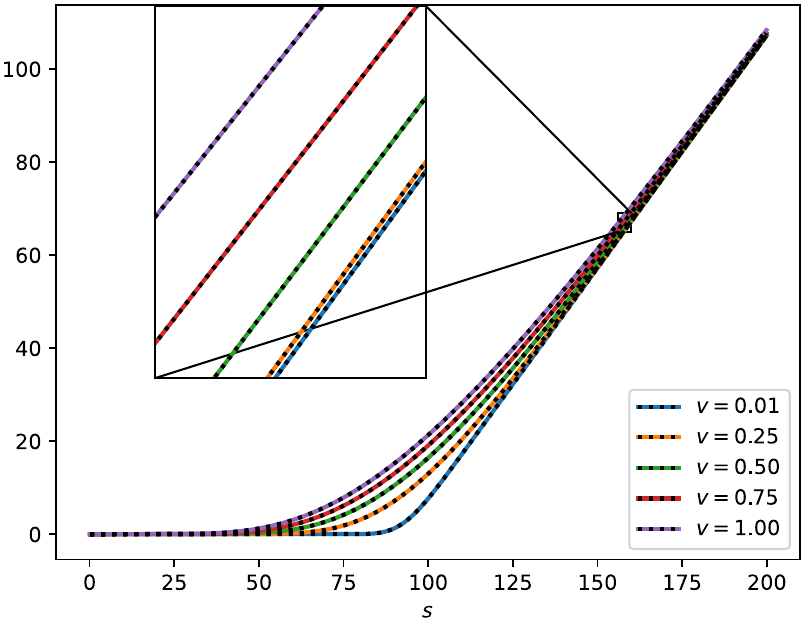}}
\caption{Cuts of prices surfaces. Numerical solution, continuous; COS solution, squares.}
\label{fig:heston-prices-cuts}
\end{figure}

In this experiment we price a vanilla call option under the Heston stochastic volatility model. The PDE model \eqref{eq:HestonPDEForward} can be written in the conservative form \eqref{SistConsCompacto} as follows:
{\small
\vspace{-0cm}
\begin{equation} \label{SistConsHeston}
    \dfrac{\partial u}{\partial t}  + \dfrac{\partial f_1}{\partial s} (u)+\dfrac{\partial f_2}{\partial v}  (u)=\dfrac{\partial g_1}{\partial s} (u_{s},u_{v})+\dfrac{\partial g_2}{\partial v} (u_{s},u_{v})+h(u),
\end{equation}
}
where the functions $f_1, f_2, g_1, g_2$ and $h$ are given by:
{\small
\vspace{-0cm}
\begin{align*}
f_1(u)&=(v-r+q)s u(s,v,t),\\
f_2(u)&=\left(\rho \sigma v - \kappa(\theta-v) + \frac12 \sigma^2 \right)u(s,v,t),\\
g_1(u_s,u_v)&=\dfrac{1}{2}s^2v u_s(s,v,t)+\rho\sigma s v u_v(s,v,t),
\end{align*}
}
{\small
\vspace{-1cm}
\begin{align*}
g_2(u_s,u_v)&=\dfrac{1}{2}\sigma^2v u_v(s,v,t),\\
h(u)&=(v-2r+q+\kappa+\rho\sigma) u(s,v,t).
\end{align*}
}
% The initial condition $u_0$, being equation \eqref{eq:hestonPayoff} for the European call option, was averaged on each volume. 
% The initial volume averages are set as
% $${\bar{u}^0}_{i,j} = \dfrac{1}{\lvert V_{ij} \rvert} \int_{V_{ij}} u_0(s)\dx{s}\dx{v}.$$
In this section we perform two vanilla call option pricing tests with the market parameters given by  $q=0.0$, $\kappa=1.5$, $\theta=0.04$, $\rho=-0.9$, $T=0.25$, $K=100$
and $\sigma=0.3, r=0.025$ for the Test 3, and  $\sigma=0.025, r=0.3$  for Test 4.
% The first set of parameters, labelled as Test 3, is taken from \cite{Foulon10}. 
The second group of parameters, denoted as Test 4, is a variation of Test 3, swapping the interest rate $r$ and the volatility of the volatility $\sigma$. In both tests, in order to minimize the numerical errors coming from the artificial boundary conditions, the computational space domain was considered as $(s,v)\in[0,800]\times[0,4]$.
% \begin{table}[!htb]
% \begin{center}
% \caption{Market data for vanilla call options under Heston model.}
% \begin{tabular}{c|c|c|c|c|c|c|c|c}
%   &$\sigma$ & $r$ & $q$  & $\kappa$ & $\theta$  & $\rho$ & $T$& $K$\\
% \hline
% Test 3& $0.3$ & $0.025$  & \multirow{2}{*}{$0.0$} & \multirow{2}{*}{$1.5$} & \multirow{2}{*}{$0.04$} & \multirow{2}{*}{$-0.9$}& \multirow{2}{*}{$0.25$} & \multirow{2}{*}{$100$}\\
% Test 4 &$0.025$ & $0.3$  &  &  &  &  &  & \\
% \end{tabular} 
% \label{tb-Heston}
% \end{center}
% \end{table}
% Figure \ref{fig:heston-prices} displays the finite volume IMEX Runge-Kutta approximations of the exact option price functions $u$ at $t=T$ corresponding to the two cases of Table \ref{tb-Heston}. 
Sections of the computed price surfaces, with planes parallel to $v=0$, are shown in Figure \ref{fig:heston-prices-cuts} in combination with the semi-analytical COS solutions. In Table \ref{tb-Heston66}, the $L_1$ errors and the $L_1$ orders of the explicit and IMEX methods for the Test 4 are displayed, respectively. Again, both numerical schemes achieve second-order accuracy in the $L_1$ norm. IMEX is able to converge using much larger time steps than the explicit method, thus it consumes much less computing time.

% Additionally, Figures \ref{fig:heston-deltas} and \ref{fig:heston-gammas} present the numerical Greeks of the numerical approximation, deltas and gammas, respectively. 
% Once again, the here proposed numerical scheme provides accurate and non-oscillatory approximations of the Greeks, even at regions of discontinuities and non-smoothness in the initial condition.

% \begin{figure}[!htb]
% \centering
% \subfigure[Test 1] {\includegraphics[height=4cm]{./images/Heston/precio_heston_65}}
% \subfigure[Test 2] {\includegraphics[height=4cm]{./images/Heston/precio_heston_convectivo_Test4}}
% \caption{Vanilla call option prices for Test 3 (left) and Test 4 (right) at $t=T$.}
% \label{fig:heston-prices}
% \end{figure}

\begin{table}[!h]
% \begin{center}
{\scriptsize
\vspace{-1cm}
\begin{tabular}{|c||c|c|c|c|}
\hline
& \multicolumn{4}{c|}{IMEX} \\
\hline
$N_1\times N_2$  & $L_1$ \text{ error} &  \text{ Order}& $\Delta t$ & Time (s)  \\
\hline
$25\times 25$ &  $8.9656\times 10^{1}$ & $--$ & $3.87\times 10^{-3}$ & $1.5\times 10^{-2}$ \\
$50\times 50$ &  $2.4203\times 10^{1}$ & $1.89$& $1.94\times 10^{-3}$ & $9.0\times 10^{-2}$\\
$100\times 100$ &  $9.3022\times 10^{0}$ & $1.38$&  $9.69\times 10^{-4}$ & $7.0\times 10^{-1}$ \\
$200\times 200$ &  $1.2578\times 10^{0}$ & $2.89$ & $4.84\times 10^{-4}$ & $6.9\times 10^{0}$ \\
$400\times 400$ &  $2.9883\times 10^{-1}$ & $2.07$ & $2.42\times 10^{-4}$  & $6.8\times 10^{1}$ \\
$800\times 800$ &  $5.9846\times 10^{-2}$ & $2.32$ & $1.21\times 10^{-4}$ & $6.3\times 10^{2}$ \\
\hline
\end{tabular} 
}
{\scriptsize
\begin{tabular}{|c|c|c|c|}
    \hline
& \multicolumn{3}{c|}{Explicit} \\
\hline
 $L_1$ \text{ error} &  \text{ Order}& $\Delta t$ & Time (s)  \\
\hline
 $9.0576\times 10^{1}$ & $--$  & $1.99\times 10^{-4}$ & $1.8\times 10^{-1}$\\
 $2.4440\times 10^{1}$ &$1.89$  & $4.99\times 10^{-5}$ & $2.7\times 10^{0}$\\
  $9.3648\times 10^{0}$ &$1.38$ & $1.25\times 10^{-5}$ & $4.4\times 10^{2}$ \\
 $1.2671\times 10^{0}$ &$2.89$ & $3.12\times 10^{-6}$ & $7.1\times 10^{2}$ \\
 $3.0089\times 10^{-1}$ &$2.07$ & $7.79\times 10^{-7}$ & $1.2\times 10^{4}$ \\
 $6.0411\times 10^{-2}$ &$2.32$  &  $1.95\times 10^{-7}$ & $1.9\times 10^{5}$ \\
\hline
\end{tabular} 
}
% \end{center}
\caption{$L_1$ errors and orders of convergence of the IMEX and explicit methods,Test 4.}
\label{tb-Heston66} 
\end{table}

\vspace{-1.5cm}
\section{Conclusions} \label{conclusions}

We have shown that finite volume IMEX Runge-Kutta numerical schemes are suitable for solving convection-diffusion PDE option pricing problems. This opens the door to the application of these schemes to numerous models in finance, even those giving rise to non-linear PDEs.
The obtained numerical schemes are highly efficient. On the one hand, large time steps can be used, avoiding the need to use small time steps enforced by the diffusion stability condition that appears when explicit schemes are considered. On the other hand, the schemes are second order accurate. 
This fact is of paramount importance in order to obtain accurate approximations of the Greeks without oscillations. We have shown that second order accuracy is preserved even when non smooth initial conditions (payoffs) are considered, which is the usual situation in finance. Thus, additional smoothing techniques for the initial condition do not need to be taken into account. Moreover, the here developed option price calculators can be extended to build numerical solvers with higher order.
% Last, but not least, in this work we provide an alternative way to compute very accurate approximations of the prices of arithmetic basket call options by means of a highly accurate and efficient multidimensional COS Fourier method. 
% These Fourier semi-analytical solutions can be also very valuable benchmarks for a broader audience that work in the development of high order schemes in the general parabolic setting, even outside the financial world.

\section*{Acknowledgements}

 % been funded by FEDER and the Spanish Government  through the coordinated Research project RTI2018-096064-B-C1; by the Junta de Andaluc\'ia research projects P18-RT-3163 and the Junta de Andalucia-FEDER-University of M\'alaga research project UMA18-FEDERJA-16; and the University of M\'alaga. 

M. Castro has has been partially supported by the grant PDC2022-133663-C21 funded by MCIN/AEI/10.13039/501100011033 and “European Union NextGenerationEU/PRTR"  and  the grant \mbox{PID2022-137637NB-C21} funded by MCIN/AEI/\-10.13039/50110001103 and “ERDF A way of making Europe”.
The other authors' research has been funded by the Spanish MINECO under research project number PDI2019-108584RB-I00 and by the grant ED431G 2019/01 of CITIC, funded by Conseller\'ia de Educaci\'on, Universidade e Formaci\'on Profesional of Xunta de Galicia and FEDER.

\end{document}